\documentclass[11pt]{article}

\usepackage[english]{babel}                        % language
\usepackage{geometry}
\usepackage[T1]{fontenc}			   % pour les mots en francais

\usepackage{theorem}{\theorembodyfont{\slshape}
  \newtheorem{theorem}{Theorem}[section]           % theorem
  \newtheorem{proposition}[theorem]{Proposition}   % proposition
  \newtheorem{corollary}[theorem]{Corollary}       % corollary
  \newtheorem{lemma}[theorem]{Lemma}               % lemma
}{\theorembodyfont{\upshape}
  \newtheorem{remark}[theorem]{Remark}             % remark
}%proofs shall not be in a 'theorem' style
\def\Proof{{\smallskip\noindent{\em Proof. }}}     % proof begin
\def\endProof{{\hfill$\Box$\medskip\noindent}}     % proof end
\newcommand{\eg}{\textit{e.g.},~}                  % latin e.g.
\newcommand{\ie}{\textit{i.e.}~}                   % latin i.e.
\newcommand{\multindice}[2]%                       % 2 lines stack
	{{\substack{{#1}\\{#2}}}}
\makeatletter                                      % underline in math mode
\newcommand*{\ubar}{\mathpalette\@ubar}
\newcommand*{\@ubar}[2]{\hbox{\b{$#1#2\m@th$}}}
\makeatother

\usepackage{amsmath,amsfonts,amssymb}              % AMS packages

\newcommand\R{{\mathbb R}}                         % real numbers
\newcommand{\e}[1]{\mathbb{E}_{#1}}                % canonical basis of R^d
\newcommand{\norme}[2][]{\left\|#2\right\|_{#1}}   % norm
\let\div=\relax\DeclareMathOperator{\div}{div}     % divergence operator
\renewcommand\P{{\mathbb P}}                       % Leray-Hopf projector
\newcommand{\NS}[1][ ]{\eqref{NS}#1}               % Navier-Stokes system
\newcommand{\cc}{K}                                % matrix \iint(u_h u_k)
\DeclareMathOperator{\id}{Id}                      % identity matrix
\DeclareMathOperator{\Tr}{Tr}                      % trace of a matrix
\newcommand{\Sphere}{\mathbb{S}}		   % unit sphere

\title{\sc New Asymptotic Profiles of Nonstationnary\\ Solutions
 of the Navier--Stokes System}

\author{
\begin{minipage}[t]{38ex}
\begin{center}Lorenzo Brandolese\end{center}\vspace*{-1.7ex}
{\sf\small%
Universit\'e de Lyon~; Universit\'e Lyon~1~;
CNRS, UMR 5208 Institut Camille Jordan,
21 avenue Claude Bernard,  {F69622} Villeurbanne Cedex, France.
\vspace*{-1.7ex}\begin{center}{\tt brandolese@math.univ-lyon1.fr}\end{center}}
\end{minipage}\qquad\begin{minipage}[t]{38ex}
\begin{center}François Vigneron\end{center}\vspace*{-1.7ex}
{\sf\small%
Centre de Math\'ematiques Laurent Schwartz, \'Ecole polytechnique,
UMR 7640 du CNRS, {F91128} Palaiseau Cedex, France.
\vspace*{-1.7ex}\begin{center}{\tt francois.vigneron@normalesup.org}\end{center}}
\end{minipage}
}

\date{To appear in Journal de Mathématiques Pures et Appliquées}

\begin{document}

\maketitle

\begin{abstract}
We show that solutions $u(x,t)$ of the non-stationnary incompressible
Navier--Stokes system in $\R^d$ ($d\geq2$) starting from mild
decaying  data $a$  behave as $|x|\to\infty$ as
a potential field:
\begin{equation}\tag{i}
u(x,t)=e^{t\Delta}a(x) + \gamma_d \, \nabla_x\left(\sum_{h,k} 
\frac{\delta_{h,k}\, |x|^2 - d x_h x_k}{d|x|^{d+2}}\,K_{h,k}(t)
\right)+\mathfrak{o}\left(\frac{1}{|x|^{d+1}}\right)
\end{equation}
where $\gamma_d$ is a constant and $K_{h,k}=\int_0^t(u_h\vert u_k)_{L^2}$
is the energy matrix of the flow.

We deduce that, for well localized data, and
for small~$t$ and large enough~$|x|$, 
\begin{equation}\tag{ii}
c\,t\,|x|^{-(d+1)}\le |u(x,t)|\le c'\,t\,|x|^{-(d+1)},
\end{equation}
where  the lower bound holds on the complementary of a set of directions,
of arbitrary small measure on $\mathbb{S}^{d-1}$.
We also obtain new lower bounds for the large time decay 
of the weighted-$L^p$ norms, extending
previous results of Schonbek, Miyakawa, Bae and Jin.
\end{abstract}

\bigskip

\renewcommand{\abstractname}{Nouveaux profils asymptotiques\\
des solutions non-stationnaires
de Navier-Stokes}
\begin{abstract}
On montre que la solution $u(x,t)$ de l'équation de
Navier--Stokes incompressible dans~$\R^d$ ($d\geq2$) issue d'une
donnée de Cauchy générique et modérément décroissante~$a$ se comporte,
pour $|x|\to\infty$, comme un écoulement potentiel
donné par la formule (i) ;
$\gamma_d$ est une constante et $K_{h,k}=\int_0^t(u_h\vert u_k)_{L^2}$
est la matrice d'énergie de l'écoulement.

On en d\'eduit que, si la donnée est bien localisée, le champ de vitesse
vérifie (ii) pour $t$ suffisament petit et $|x|$ assez grand.
La borne inférieure est valable sur le complémentaire d'un ensemble
de directions, de mesure arbitrairement petite dans $\mathbb{S}^{d-1}$.
On obtient aussi de nouvelles bornes inférieures du taux de décroissance en
temps grand des moments de la solution dans $L^p$ qui 
étendent des résultats antérieurs
de Schonbek, Miyakawa, Bae and Jin.
\end{abstract}

\vspace{1ex}
\begin{flushleft}
Lower bound estimates. Mild solutions to the Navier-Stokes system.
Asymptotic separation of variables.
 Peetre weight.
\end{flushleft}

% ----------------------------------------------------------------------------------

\section{Introduction}
Let $a$ be a divergence-free vector field in $\R^d$ ($d\geq2$).
We consider the Cauchy problem for the Navier--Stokes equations~:
\begin{equation}\label{NS}\tag{NS}
\begin{cases}
\partial_t u -\Delta u +(u\cdot\nabla)u=-\nabla p,\\
\div u = 0,\\
u(x,0)=a(x).
\end{cases}
\end{equation}
The unknowns are the velocity field~$u=(u_1,\ldots,u_d)$ and the
pressure $p$. The problem has to be solved on $\R^d\times[0,+\infty)$
or at least on $\R^d\times[0,T)$ for some $T>0$.

\bigskip
Because of their parabolic nature, the Navier--Stokes equations feature
an infinite-speed propagation effect in the space variable.
This phenomenon is usually described by the fact that compactly supported
initial data give rise to solutions which immediately have non-compact support.
On the other hand, because of the pressure, which can be eliminated from
the equations only by applying a non-local operator,
the solutions of the Navier--Stokes equations have quite a different behavior
as $|x|\to\infty$ from that of solutions of non-linear heat equations.

The main purpose of this paper is to study such asymptotic behavior.
For example, we address the following problem.
Assume that, at the beginning of the evolution, the fluid is at rest outside
a bounded region (say, $a\in C^\infty_\sigma(\R^d)$, the space of smooth, solenoidal
and compactly supported vector fields).
\emph{At which velocity will the fluid particles that are situated far from that
region start to move~?}

We will obtain sharp answers to this and related questions
by constructing new asymptotic profiles 
of solutions to (NS), predicting
the pointwise behavior of $u$ as $|x|\to+\infty$.

A few asymptotic profiles of solutions to the Navier--Stokes equations
in the whole space are, in fact, already known.
For example, 
F.~Planchon \cite{Pla98}, studied self-similar profiles.
However, his results cannot be used
in the case of initial data decaying at infinity faster than $|x|^{-1}$,
since the only possible self-similar profile, in this case,  would be the zero function.
For faster decaying data,
the asymptotic profiles of 
 A.~Carpio \cite{Car96}, Y.~Fujigaki, T.~Miyakawa~\cite{FujM01},
Miyakawa, Schonbek \cite{MS}, 
T.~Gallay, E.~Wayne \cite{GW1} and Cannone, He, Karch~\cite{CHK06}
provide valuable information
about the large-time behavior of the velocity field.
However, in all these works the asymptotics are obtained by computing
some spatial norms of expressions involving the solution.
The limitation of this approach is that most of the information on the
pointwise behavior of the velocity field is lost.

Our method is different, and consists in proving that, asymptotically,
the flow behaves as a linear combination of functions of separate
variables $x$ and $t$.

\bigskip

Our profiles imply that, without external forces, the flow
associated with decaying initial data behaves at infinity as a potential
field, with a generalized Bernouilli formula relating the pressure to the
energy matrix $(u_h\vert u_k)_{L^2}$ of the flow.
This illustrates the fact that the spatial behavior at infinity of the flow
is almost time-independent, contrary to the temporal asymptotic, which is
known to be influenced by spatial decay.

\subsection*{Notations}

\begin{enumerate}
\item
We denote by $L^\infty_\vartheta$ the space of all measurable functions
(or vector fields) $f$ on $\R^d$, such that~:
$$\|f\|_{L^\infty_\vartheta}=
\underset{x\in\R^d}{\operatorname{ess~sup}}\:(1+|x|)^{\vartheta}\,|f(x)|<+\infty.$$
The space ${\cal C}_w\left([0,T);L^\infty_\vartheta\right)$ is made of functions
$u(x,t)$ such that $u(t)\in L^\infty_\vartheta$ for all $t\geq0$ and
$$\begin{cases}
\lim\limits_{t'\to t}\|u(t')-u(t)\|_{L^\infty_\theta}=0 &
\text{if}\enspace t>0,\\[1ex]
u(t) \underset{t\to0}{\rightharpoonup} u(0) &
\text{in the distributional sense}.\end{cases}$$
\item
For positive $\phi$,
the notation $f(x,t)={\cal O}_t\left(\phi(x)^{-1}\right)$ means that
$|\phi(x)f(x,t)|\le C_t$,  for some function  $t\mapsto C_t$, possibly 
growing as $t\to\infty$, but locally bounded.
\item
The solution of the heat equation is
$$e^{t\Delta}a(x) =
(4\pi t)^{-d/2} \int_{\R^d} e^{-\frac{|x-y|^2}{4t}}\,a(y)\,dy.$$
\item
We also adopt the standard Kronecker symbol:  $\delta_{i,j}=1$ if $i=j$,
and $\delta_{i,j}=0$ otherwise.
\end{enumerate}

\bigskip
Our starting point is the following well known result
(see \cite{Lem02}, chapter 25).
\begin{theorem} 
\label{theorem-intro}
Let $d\ge2$. There exists a constant $\gamma>0$ such that for any
divergence-free vector field $a\in L^\infty(\R^d)$, one can find
$$T\ge \gamma\min\left\{1;\|a\|^{-2}_{L^\infty}\right\}$$
and a unique mild solution $u\in {\cal C}_w\left([0,T];L^\infty\right)$ of \NS[].
This solution $u$ is smooth for~$t>0$.
Moreover, if $a$ belongs to $L^\infty_\vartheta$
for some $\vartheta\geq0$, then we also have~:
$$u\in {\cal C}_w\left([0,T];L^\infty_{\ubar\vartheta}\right)$$
with $\ubar\vartheta=\min\{\vartheta\,;\,d+1\}$.
\end{theorem}

This conclusion can be restated in a slighlty different
way (see also \cite[Proposition~3]{Vig05})~:
$$u(x,t)=e^{t\Delta}a + {\cal O}_t\left((1+|x|)^{-\min\{2\vartheta\,;\, d+1\}}\right)$$
on $[0,T]\times\R^d$.

\subsection*{Asymptotic behavior of local solutions}
We can now state our first main result.
Let us introduce the energy matrices~:
\begin{equation}
\label{enem}
{\cal E}_{h,k}(t)=\int_{\R^d}(u_hu_k)(y,t)\,dy
\qquad\text{and}\qquad
\cc_{h,k}(t)=\int_0^t\!\!\int_{\R^d}(u_hu_k)(y,s)\,dy\,ds.
\end{equation}
The following theorem describes the asymptotic profile of local solutions.

\begin{theorem}
\label{theorem1}
For $\vartheta>\frac{d+1}{2}$ and an initial datum $a\in L^\infty_\vartheta$,
let $u\in {\cal C}_w\big([0,T];L^\infty_{\ubar\vartheta}\big)$ be the solution of \NS given
by the preceding theorem.
The following profile  holds for $|x|\to+\infty$~:
\begin{equation}
 \label{profile}
u(x,t)=e^{t\Delta}a(x)+ \nabla\Pi(x,t)
+{\cal O}_t\left(|x|^{-\min\{2\vartheta \,;\, d+2\}}\right)
\end{equation}
where $\Pi(x,t)$ is given by~:
\begin{equation}
\Pi(x,t)= \gamma_d\sum_{h,k}\left(\frac{\delta_{h,k}}{d\,|x|^d}
-\frac{x_h x_k}{|x|^{d+2}}\right)\cdot\cc_{h,k}(t)
\end{equation}
and $\gamma_d= \pi^{-d/2} \, \Gamma(\frac{d+2}{2})$.
If, moreover, the first and second order derivatives of~$a$
belong to~$L^\infty_\vartheta$, then there exists a constant~$p_0$ such that
the following profile holds for $t>0$~:
\begin{equation}
\label{profile pressure}
p(x,t)=p_0-\gamma_d\sum_{h,k}\left(\frac{\delta_{h,k}}{d\,|x|^d}
-\frac{x_h x_k}{|x|^{d+2}}\right)\cdot {\cal E}_{h,k}(t)
+{\cal O}_t\left(|x|^{-\min\{2\vartheta-1 \,;\, d+1\}}\right)
\end{equation}
\end{theorem}

\begin{remark}
This theorem essentially says that, for \emph{mild decaying\/} data
(this is the meaning of the assumption $a\in L^\infty_\vartheta$, with $\vartheta>\frac{d+1}{2}$),
$$ u(x,t)\sim e^{t\Delta}a(x)+\nabla\Pi(x,t), \qquad
    \hbox{as $|x|\to\infty$}.$$
In other words, the
flow behaves at infinity as the solution
of the heat equation plus a potential field at infinity.
In particular, if follows that for \emph{fast decaying\/} data ({\it i.e.\/} when
$\vartheta>d+1$, we simply have
$$u(x,t) \sim \nabla\Pi(x,t),\qquad  \hbox{as $|x|\to\infty$},$$
since the linear evolution can be included in the lower order terms.

Theorem~\ref{theorem1} does not cover the case of \emph{slowly decaying\/} data
({\it i.e.\/}, the case $\vartheta\le \frac{d+1}{2}$). The spatial 
asymptotic of those slowly decaying solutions (including self-similar solutions)
has a different structure and cannot be constructed with the same method.
We should consider it in an independent paper.
\end{remark}

\begin{remark}
The decay of the remainder in~\eqref{profile} cannot exceed $|x|^{-d-2}$. Indeed,
\NS being invariant by translation, the choice of the origin
is arbitrary and one can easily check that $$\nabla\Pi(x-x_0,t)-
\nabla\Pi(x,t)$$ decays at infinity as $|x|^{-d-2}$ if $\Pi\not\equiv0$
and $x_0\neq0$.

Even if $u(x,t)$ develops a singularity in finite time,
the potential field in \eqref{profile} will  remain uniformly bounded
away from the origin~:
$$|\nabla\Pi(x,t)|\leq C\norme[L^2]{a} t |x|^{-d-1}.$$
However, the above result provides no  information about the singularity
itself, nor does it prevent it from appearing~:
as long as the solution is smooth, the remainder of~\eqref{profile} compensates
for the singularity at the origin of $\nabla\Pi(x,t)$.
\end{remark}

\begin{remark}
The above profile for the pressure has some analogy with Bernoulli's formula for potential flows~:
$$p=p_0+\frac{1}{2}\rho U^2.$$
Such a formula holds rigorously for the stationary Euler equation with no external
force, but this identity can also be useful when
dealing with high Reynolds flows around aerodynamical bodies 
(see, \eg
the description of the Prandtl laminarity theory in \cite[Chapter~9]{guyon-hulin}).
\end{remark}

\bigskip
The asymptotic profiles of Theorem~\ref{theorem1} are meaningful
when the leading term does not vanish identically.
It turns out that this is the case for generic solutions.
Indeed, the next result provides a necessary and sufficient condition for $\nabla\Pi$
to be identically zero.

\begin{proposition}
\label{theorem2}
Let $u$ as in Theorem~\ref{theorem1} and $\cc=(\cc_{h,k})$.
For any~$t\in [0,T]$, the homogeneous function~$x\mapsto\nabla\Pi(x,t)$
vanishes identically on $\R^d$ if and only if
the matrix $K(t)$ is proportional to the identity matrix, \ie
\begin{equation}
\label{orth}
\forall \,h,k\in\{1,\ldots,d\},\qquad
\cc_{h,k}(t)=\alpha(t)\,\delta_{h,k}
\end{equation} 
with $\alpha=\frac{1}{d}\Tr\cc$.
\end{proposition}

This shows that $\nabla\Pi$ does not vanish
for generic flows.
Conditions \eqref{orth} also occur in the paper of
T.~Miyakawa and M.~Schonbek~\cite{MS}.
It is shown therein that a high decay rate
of the energy of the flow for large time is
essentially equivalent to~\eqref{orth} holding in the limit~$t\to+\infty$.

Such orthogonality relations can also be described in terms of
vanishing moments of the vorticity $\omega=\hbox{curl}\, u$
of the flow
(see, {\it e.g.\/}  \cite{GW1}, \cite{GW2}).
Focusing on the vorticity, in fact, has crucial advantages in the study of the large
time behavior of solutions, especially in the two-timensional case.
We refer {\it e.g.\/}
to the recent work of Gallay and Wayne \cite{GW3}
on the global stability results of vortex solutions.

On the other hand, the large space behavior of the vorticity
is less interesting than that of the velocity field.
This can be shown by taking the 
$\hbox{curl}(\cdot)$ operator term-by-term in formula~\eqref{profile}:
the term $\hbox{curl}(\nabla\Pi)$ identically vanishes.
The physical interpretation of this remark is the following:
if we start with an initial datum with localized vorticity, then 
the vorticity will remain localized as far as the solution exists
(this fact, of course, was already known).

In principle, it would be possible
to extend formula~\eqref{profile}
and to write a higher-order asymptotic for $u$ as $|x|\to\infty$.
The above observation allows us to predict that all the higher-order terms
of the expansion of the velocity field must be $\hbox{curl}$-free.
Otherwise, there would be a limitation on the decay rate of $\omega$ as $|x|\to\infty$
and this would contredict {\it e.g.\/} the results of~\cite{GW2}, \cite{KuT05}.
In other words, all the higher-order terms of the expansion of $u$ should be gradients.

\subsection*{Large time asymptotics}

Under a suitable smallness assumption such as 
$$\underset{x\in\R^d}{\operatorname{ess~sup}}\: |x|\,|a(x)|\le \varepsilon_0,$$
one can take $T=+\infty$ in Theorem~\ref{theorem-intro} (see, \eg\cite{Mi00}).
Moreover, the localization property of the flow persists uniformly.
One has~:
$$|u(x,t)|\leq C \, (1+|x|)^{-\alpha}  \, (1+t)^{-\beta/2}$$
for any $\alpha,\beta\geq0$ such that $\alpha+\beta\leq\min\{\vartheta \,;\,
d+1\}$. When $\vartheta=d$ or $d+1$, one needs the additional assumption that
the above estimate already holds for $e^{t\Delta}a$ (see \cite{Bra04i}).
In particular, these estimates imply that, for large $t$~:
$$\norme[{L^2([0,t],L^2)}]{u}^2 \le
\begin{cases}
C, &\hbox{if $\vartheta>\frac{d+2}{2}$}\\
C_\epsilon\, t^{-\vartheta+\frac{d+2}{2}+\epsilon}, &\hbox{if $\vartheta\le \frac{d+2}{2}$},
\end{cases}
$$
for all $\epsilon>0$ (this bound also holds for $\epsilon=0$, $\vartheta\not=\frac{d+2}{2}$,
but we will not use this fact).

\smallskip
We can now give our asymptotic profile for global solutions.

\begin{theorem}
\label{theorem3}
Given $\vartheta>\frac{d+1}{2}$,
let $u(x,t)$ be a solution of \NS on $\R_+\times\R^d$ such that
\begin{equation}
\label{assumption_largetime}
|u(x,t)|\leq C_0 \, (1+|x|)^{-\alpha}  \, (1+t)^{-\beta/2}
\end{equation}
for any $\alpha,\beta\geq0$ such that $\alpha+\beta\leq\ubar\vartheta$.
Then,
\begin{equation}
\label{detailed profile}
u(x,t)= e^{t\Delta}a(x) + \nabla\Pi(x,t) + |x|^{-d-1}
E\!\left(\frac{x}{\sqrt{t+1}}\,;\,t\right)+ {\cal R}(x,t)
\end{equation}
with the following estimates~:
\begin{equation}
\label{double bound}
|E(x,t)|\leq C e^{-c|x|^2}  \norme[{L^2([0,t],L^2)}]{u}^2
\end{equation}
and, for any $0\le\alpha\le \min\{1,\vartheta-\frac{d+1}{2}\}$, and all $t\ge1$,
\begin{equation}
\begin{aligned}
\label{remainder_asympt}
|{\cal R}(x,t)| &\leq C_\alpha \, |x|^{-d-1-\alpha}\, t^{-\frac{1}{2}+\frac{\alpha}{2}}, 
	\qquad\hbox{\qquad if $\vartheta>\frac{d+3}{2}$},\\
|{\cal R}(x,t)| &\leq C_{\alpha,\epsilon} \, |x|^{-d-1-\alpha} \,t^{\frac{d+2+\alpha}{2}-\vartheta+\epsilon},
 \qquad\hbox{if $\frac{d+1}{2}<\vartheta\le\frac{d+3}{2}$}.
\end{aligned}
\end{equation}
\end{theorem}

Due to the form of the remainder terms, it seems impossible to
obtain a  description of the pointwise behavior of $u$ for large but fixed $|x|$, and $t\to\infty$.
Conclusion~\eqref{detailed profile} is interesting only for $(x,t)$ such that
$|x|\ge C\sqrt{t+1}$.
For those points, this profile 
provides more information 
than those in~\cite{FujM01} or \cite{CHK06}
(on the other hand, our assumptions are necessarily more stringent).
%Note that all the singularities compensate
%each other in the region $|x|\leq C\sqrt{t+1}$, which allows
%the velocity field to tend uniformly to zero for large~$t$.

\subsection*{Applications}

For smooth and fast decaying initial data, according to Theorem~\ref{theorem1}
one has 
\begin{equation}
\label{Miya-bound}
|u(x,t)|\le C_t \, (1+|x|)^{-(d+1)}.
\end{equation}
Theorem~\ref{theorem1} allows us to answer the more subtle problem
of the validity of the corresponding {\em lower bound\/} to \eqref{Miya-bound}.

\smallskip
A first difficulty is the following:
the upper bound  ensures that $u(\cdot,t)$ is integrable, so that
the divergence-free condition implies
$$\forall t>0,\quad\forall j\in\{1,\ldots,d\},\qquad
\int_{\R^d}u_j(x,t)\,dx=0.$$
In particular, since $u$ is smooth for $t>0$, no uniform lower bound by a given
positive function can hold.
However,
non-uniform and anisotropic lower bounds do hold, even if the initial data
is rapidly decreasing.

More precisely, for generic flows (i.e. if we exclude flows with special symmetries)
starting from fast decaying data,
we will prove that for some small $t_0>0$ (depending only on the initial datum), and for $j=1,\ldots,d$,
\begin{equation}
\label{lower-upper bounds}  
  c\,t\,|x|^{-(d+1)} \le  |u_j(x,t)|\le c'\,t\,|x|^{-(d+1)}, \qquad c,c'>0
\end{equation}
for all $t\in (0,t_0]$ and all $|x|\ge C/\sqrt t$, with $x$ outside a small set of exceptional directions,
along which the decay can be faster.
In other words, the constant $c$ in~\eqref{lower-upper bounds} is independent of $t$ or $|x|$,
but does depend on the direction ${x}/{|x|}$
(see Theorem~\ref{corollary1} below for a more precise statement).
For example, we will see that in dimension two the exceptional set is made of at most \emph{six directions}.
The remarkable fact is that the above lower bound holds {\it e.g.\/} for compactly supported
data (that is, even without assuming that $|a(x)|\ge  c(1+|x|)^{-d-1}$).
In particular, this allows us to improve the previously known results
(see \eg \cite{BraM02}, \cite{Lem02})
on the instantaneous spatial spreading property of highly localized flows.

\bigskip
For generic global strong solutions, Theorem~\ref{theorem3} implies 
various lower bounds.
More precisely, starting from a fast decaying initial datum, we get,
for all $0\le \alpha\le d+1$, and large~$t$~:
\begin{equation}
\label{large time lower0}
\norme[{L^\infty_\alpha}]{u(t)} \ge c\,t^{-(d+1-\alpha)/2}.
\end{equation}
This result is a converse to Miyakawa's property~\eqref{assumption_largetime}.

In the same spirit,
Theorem~\ref{theorem3} can be applied to 
estimate the decay of the moments of  the solutions~:
for all $1\le p<\infty$ and  $\alpha\ge0$ such that
\begin{equation}
\label{localization restriction}
\alpha+\frac{d}{p}<d+1,
\end{equation}
we obtain, for large~$t$,
\begin{equation}
\label{large time lower}
\norme[{L^p}]{\,(1+|x|)^\alpha\,u(t)} \geq
c\,t^{{-\frac{1}{2}(d+1-\alpha-\frac{d}{p})}}. 
\end{equation}
This  lower bound seemed to be  known only 
in a few particular cases (namely,  $p=2$ and $0\le\alpha\le2$, see \cite{Sch91}, \cite{GW2}, \cite{BaJ05},
or $1\le p\le \infty$ and $\alpha=0$, see~\cite{FujM01}).
The corresponding upper bounds to~\eqref{large time lower}, starting with the work
of M.~E. Schonbek,  have been studied by many authors.
See \cite{KuT05} for a quite general result.

In some sense, the restriction~\eqref{localization restriction} on the parameters could be removed,
since~\eqref{lower-upper bounds} implies that
for generic solutions,  one has
$$\norme[{L^p}]{\,(1+|x|)^\alpha\,u(t)}=\infty$$
whenever $\alpha+\frac{d}{p}\ge d+1$
(see also~\cite{BraM02}).

\bigskip
Our results can also be applied to the study of the anisotropic decay of the
velocity field.
In the whole space, we show that not too stringent anisotropic assumptions 
on the decay of the data will be conserved by the flow.
We also show that, if the initial data is well localized in $\R^d$,
then the  flow decays faster than $(1+|x|)^{-(d+1)}$ as soon as one
component does. This prevents  localized flows in $\R^d$
from having a really anisotropic decay.

The situation can be different in other unbounded domains.
For example, we will briefly discuss the 
case of the half plane $x_d>0$, with Neumann boundary conditions,
and show that, in this case, generic flows have a genuinely anisotropic decay.

\bigskip
\subsection*{The asymptotic separation of variables method}
The proof of~\eqref{profile} relies on a new, simple method
that is a sort of ``asymptotic separation of variables''.
We can summarize it as follows:
one starts writing the Navier--Stokes equation in the usual integral form
\begin{equation}\label{IE0}
u(t)=e^{t\Delta}a-\displaystyle\int_0^t e^{(t-s)\Delta}\,\P\div(u\otimes u)(s)\,ds,
\end{equation}
where $\P$ is the Leray-Hopf projector onto the divergence-free vector fields~:
$$\P f=f-\nabla\Delta^{-1}(\div f).$$
Then we use a classical decomposition
of the nonlinear term (see {\it e.g.\/},  \cite{Sch91}, \cite{BraM02})
$$ (u\otimes u)(x,t)=\biggl(\int_{\R^d}(u\otimes u)(y,t)\,dy\biggr) g(x) +
v(x,t),$$
where $g$ denotes the standard gaussian function, and $v$ is defined through this formula.

Since $\displaystyle\int_{\R^d}v(x,t)\,dx=0$, the function
$e^{t\Delta}\,\P\div v$ behaves at infinity 
better than
the previous non-linearity $e^{t\Delta}\P\,\div(u\otimes u)$~:
its contribution can be included in the remainder terms.
The next step consists in
observing that
the kernel of $e^{t\Delta}\P\,\div$ behaves, as $|x|\to\infty$,
as a time independent homogeneous tensor $H(x)$.
Then we show that by applying $e^{(t-s)\Delta}\P\,\div$ to 
a matrix of the form ${\cal E}(s)g$, where the coefficients
of ${\cal E}(s)$ depend only on time,
we get $H(x)\cdot{\cal E}(s)$, plus some lower order terms.
A time integration then yields a principal part for the velocity field
of the form $H(x)\cdot K(t)$, as $|x|\to\infty$.
An explicit computation of this product
provides the expression for $\Pi(x,t)$ in \eqref{profile}.

We point out that the above strategy is not specific to the Navier--Stokes
equations, but can be adapted to obtain the spatial asymptotics
for more general models. What one essentially needs for its application
are sufficiently explicit expressions (or sharp estimates) for the kernels
of the operators involved.

\subsection*{Structure of the article}

Our main results are Theorem~\ref{theorem1}, its companion Theorem~\ref{theorem3},
and~Theorem~\ref{corollary1}.
Corollary~\ref{corollary4} also has some interest, since it extends a few
results in the existing literature and its proof is very short.
This paper is organized as follows.
The proof of Theorem~\ref{theorem1} and
Theorem~\ref{theorem3} is contained in Sections~\ref{sub2.4}--\ref{sub2.5},
after we have prepared some preliminary estimates.
In Section~\ref{sub2.6} we establish Proposition~\ref{theorem2},
in a slightly more complete form.
The remaining part of the paper is devoted to applications~:
in Section~\ref{sub3.1} we give a precise statement and a proof of~\eqref{lower-upper bounds}.
Section~\ref{sub3.2} contains the proof of~\eqref{large time lower0} and \eqref{large time lower}.
The last sections deal with the anisotropic decay of solutions.

\section{Proof of the main results}
\label{section2}

Let us now focus on the proof of the above results.

\bigskip
We shall use the following notations for the kernel of the
convolution operator $e^{t\Delta}\P\div$ ~:
$$F_{j;h,k}(x,t) = \int_{\R^d} i e^{-t|\xi|^2+ix\cdot\xi}
\left(\frac{1}{2}\left[\xi_h \delta_{j,k}+\xi_k\delta_{j,h}\right]
-\frac{\xi_j\xi_h\xi_k}{|\xi|^2}\right) \, \frac{d\xi}{(2\pi)^d}\cdotp$$
According to \eqref{IE0}, the $j^\text{th}$~component of \NS can
therefore be written as
\begin{equation}
\label{IE}
 u_j(t)=e^{t\Delta}a_j-
\sum_{h=1}^d \sum_{k=1}^d\int_0^t F_{j;h,k}(t-s)*(u_{h}u_{k})(s)\,ds.
\end{equation}

\medskip

This kernel is related to the standard gaussian function $g(x)=(4\pi)^{-d/2}e^{-|x|^2/4}$
in the following way. One has~:
$$F_{j;h,k}(x,t)=F^{(1)}_{j;h,k}(x,t)+F^{(2)}_{j;h,k}(x,t)$$
with
$$F^{(1)}_{j;h,k}(x,t)=\frac{1}{2}\left[(\partial_h
g_t)\delta_{j,k} + (\partial_k g_t)\delta_{j,h} \right],
\qquad
F^{(2)}_{j;h,k}(x,t)=\int_t^\infty \partial_j\partial_h\partial_k
g_s(x)\,ds$$
Note that $F_{j;h,k}^{(1)}=F^{(1)}_{j;k,h}$ and  $F_{j;h,k}^{(2)}=F^{(2)}_{j;k,h}$
accordingly to the fact that only the symmetrical kernel has a physical meaning;
$g_t(x)=t^{-d/2}g(x/\sqrt t)$ is the fundamental solution of the heat
equation.

\subsection{Some elementary computations on $F$}
\label{sub2.1}

We shall need time-independent asymptotics of $F$,
valid in the region where $|x|^2\gg t$.
\begin{lemma}
\begin{subequations}
There exist two positive constants $C$ and $c$ that depend only on $d$,
and a family of smooth functions $\Psi_{j;h,k}$ satisfying
\begin{equation}
\label{bound-Psi}
 |\Psi_{j;h,k}(x)|+|\nabla\Psi_{j;h,k}(x)|\le Ce^{-c|x|^2}
\end{equation}
such that~:
\begin{equation}\label{F}
F_{j;h,k}(x,t) = \gamma_d\left(\frac{\sigma_{j,h,k}(x)}{|x|^{d+2}}
-(d+2)\frac{x_jx_hx_k}{|x|^{d+4}}\right)
+|x|^{-(d+1)}\Psi_{j;h,k}\biggl(\frac{x}{\sqrt{t}}\biggr)
\end{equation}
with $\gamma_d= \pi^{-d/2} \, \Gamma(\frac{d+2}{2})$ and
$\sigma_{j,h,k}(x)=\delta_{j,h}x_k+\delta_{j,k}x_h+\delta_{h,k}x_j$.
\end{subequations}
\end{lemma}
\begin{remark}
Note that $x=0$ is not a singular value of $F$ ; indeed,
$F_{j;h,k}$ is a ${\cal C}^\infty$ function on $\R^d\times]0;+\infty[$
and one may immediately check on the Fourier transform that $F_{j;h,k}(0,t)=0$.
Actually, for $|x|^2\leq t$, the following computations also imply
that~:
$$|F_{j;h,k}(x)| \leq
C\, \frac{\delta_{j,h}|x_k|+\delta_{j,k}|x_h|+\delta_{h,k}|x_j|}{t^{(d+2)/2}}
+\mathcal{O}\left(\frac{|x|^2}{t^{(d+3)/2}}\right).$$
\end{remark}
\begin{remark}\label{notations_noconflict}
In Theorem \ref{theorem3}, one has~:
$$ E_j(y,t) = - \sum_{h,k} \int_0^t\!\!\int_{\R^d}
(u_hu_k)(z,t-s)\,\Psi_{j;h,k}\!\left(\sqrt{\frac{t+1}{s+1}}\:y\right)
\,dz\,ds$$
with the $\Psi_{j;h,k}$ given by this lemma.
\end{remark}

\Proof
For all indices $j,h,k$ (distinct or not) in~$\{1,\ldots,d\}$, one has~:
$$F^{(1)}_{j;h,k}(x,t)=-\frac{\delta_{j,k}\,x_h+\delta_{j,h}\,x_k}
{4(4\pi)^{d/2}\,t^{(d+2)/2} }e^{-|x|^2/4t},$$
thus $F^{(1)}_{j;h,k}(x,t) =
|x|^{-(d+1)}\,\Psi^{(1)}_{j;h,k}\left(x/\sqrt{t}\right)$,
with
$$\Psi^{(1)}_{j;h,k}(x)=-
2^{-d-1}\pi^{-d/2}(\delta_{j,k}\,x_h+\delta_{j,h}\,x_k)|x|^{d+1}e^{-|x|^2/4}.$$

Let us introduce
$\sigma_{j,h,k}(x)=\delta_{j,h}x_k+\delta_{j,k}x_h+\delta_{h,k}x_j$.
One also has~:
$$F^{(2)}_{j;h,k}(x,t)=\int_t^\infty \left(
\frac{\sigma_{j,h,k}(x)}{(2s)^2}-\frac{x_jx_hx_k}{(2s)^3}\right) g_s(x) \,ds.$$
The change of variable $\lambda=|x|/\sqrt{4s}$ gives
$g_s(x)=\pi^{-d/2}|x|^{-d}\lambda^d e^{-\lambda^2}$, and therefore~:
$$F^{(2)}_{j;h,k}(x,t)= 2\pi^{-d/2}\int_0^{|x|/\sqrt{4t}} \left(
\frac{\sigma_{j,h,k}(x)}{|x|^{d+2}}\lambda^{d+1} -
\frac{2x_jx_hx_k}{|x|^{d+4}}\lambda^{d+3}
\right) e^{-\lambda^2} \,d\lambda.$$

The following formula provides information when $A=|x|/\sqrt{t}\gg1$~:
$$\int_0^A \lambda^{d+n}e^{-\lambda^2}d\lambda =
\frac{1}{2}\Gamma\left(\frac{d+n+1}{2}\right)-
\int_A^\infty \lambda^{d+n}e^{-\lambda^2}  d\lambda.$$
This leads to~:
$$\pi^{d/2}F^{(2)}_{j;h,k}(x,t)=
\frac{\sigma_{j,h,k}(x)}{|x|^{d+2}}\,\Gamma\left(\frac{d+2}{2}\right)
-\frac{2x_jx_hx_k}{|x|^{d+4}}\,\Gamma\left(\frac{d+4}{2}\right)
+|x|^{-(d+1)}\,\Psi^{(2)}_{j;h,k}\left(\frac{x}{\sqrt{t}}\right)$$
with $$
\Psi^{(2)}_{j;h,k}(x)=-\frac{2\sigma_{j,h,k}(x)}{|x|}
\int_{|x|/2}^\infty \lambda^{d+1}e^{-\lambda^2}  d\lambda
+\frac{4x_jx_hx_k}{|x|^{3}}
\int_{|x|/2}^\infty \lambda^{d+3}e^{-\lambda^2}  d\lambda.$$
Conclusion \eqref{F} follows immediately from the well known
formula $\Gamma(z+1)=z\Gamma(z)$.
The bounds on $\Psi_{j;h,k}$ and its derivatives are also obvious. 

\endProof

\bigskip
The second valuable property of $F$ is that the convolution with the standard
gaussian function is equivalent to a shift in time.
\begin{lemma}\label{gshift}
For all $t>0$ and $x\in\R^d$, one has~:
\begin{equation}
\bigl(F_{j;h,k}(\cdot,t)*g\bigr)(x)=F_{j;h,k}(x,t+1).
\end{equation}
\end{lemma}

\Proof
Since $(g_t)_{t\geq0}$ is a convolution semi-group, \ie
$g_t\ast g = g_{t+1}$, one has~:
$$F_{j;h,k}^{(1)}(\cdot,t)\ast g = \frac{1}{2}\left[(\partial_h
g_t)\delta_{j,k} + (\partial_k g_t)\delta_{j,h} \right] \ast g
= F_{j;h,k}^{(1)}(\cdot,t+1)$$
and $\displaystyle
F^{(2)}_{j;h,k}(\cdot,t)\ast g=\int_t^\infty \partial_j\partial_h\partial_k
g_{s+1}(x)\,ds=\int_{t+1}^\infty \partial_j\partial_h\partial_k
g_{s}(x)\,ds=F^{(2)}_{j;h,k}(\cdot,t+1)$.

\endProof

\bigskip
Let us finally recall a classical estimate of the $L^1$ norm of the kernel.
\begin{lemma}
\label{lemma classical}
There exists a constant $C>0$ such that
\begin{equation}\label{L1F}
\forall t>0,\qquad \norme[L^1]{F(\cdot,t)} \leq C \, t^{-1/2}.
\end{equation}
\end{lemma}
\Proof
This follows from~\eqref{F}.

\endProof

\subsection{Decomposition of the non-linear term}
\label{sub2.2}
Theorems~\ref{theorem1} and \ref{theorem3} rely on 
a suitable decomposition of the non-linear term.
A similar decomposition has been previously used by
M.~Schonbek~\cite{Sch91} to prove lower bounds on the large-time
decay of the $L^2$-norm of the flow.
This part of the computations is common to both proofs.

\medskip
Let us first explain the decomposition on a gaussian non-linearity.
If $g$ denotes the standard gaussian function and $g_t$ the fundamental
solution of the heat equation, one sets~:
$$g_t^2(x) = \left(\int_{\R^d}g_t(y)^2 dy\right) g + \triangle(x,t).$$
The remainder $\triangle(x,t)$ has a mean value of zero~:
$$\int_{\R^d}\triangle(x,t)dx=0.$$
For fixed $x\in\R^d$, this approximation scheme behaves badly if $t\to0$
or $t\to+\infty$ ;  indeed, a simple computation leads to
$$\triangle(x,t)=\left\{1-\left(\frac{t}{2}\right)^{d/2}
e^{(2-t)|x|^2/4t}\right\}g_t^2(x).$$
But this computation also shows that such approximation scheme
is satisfactory at least when
$$t\simeq2 \;\enspace\text{or}\enspace\;
|x|^2\simeq2d\,\frac{t}{t-2}\ln\frac{t}{2}\, ,$$
{\it i.e.\/} when $\triangle(x,t)$ is close to zero.
We now perform a similar decomposition for the non-linearity
in \NS[.]

\bigskip
Now let $a$ be a $L^\infty_\vartheta$ divergence-free vector field and
$u\in {\cal C}_w\big([0,T];L^\infty_{\ubar\vartheta}\big)$ be
the solution of \NS given by Theorem~\ref{theorem-intro}, starting from
these initial data.
Recall that $\ubar\vartheta=\min\{\vartheta;d+1\}$.
Since $\vartheta>\frac{d}{2}$, one has $L^\infty_{\ubar\vartheta}\subset L^2$
and hence the energy matrix
$${\cal E}_{h,k}(t)=\int_{\R^d}(u_hu_k)(y,t)\,dy$$
is well defined.
Consistently with the preceding approximation scheme, let us define $v_{h,k}$ by
\begin{equation}
(u_hu_k)(x,t)={\cal E}_{h,k}(t) g(x) + v_{h,k}(x,t).
\end{equation}
Thanks to Lemma~\ref{gshift}, the integral equation \eqref{IE}
is hence equivalent to~:
$$ u_j(t)=e^{t\Delta}a_j-
\sum_{h,k} \int_0^t {\cal E}_{h,k}(s) \: F_{j;h,k}(t+1-s)\,ds
-\sum_{h,k}\int_0^t v_{h,k}(s)\ast F_{j;h,k}(t-s)\,ds.$$
The time-independent asymptotic \eqref{F} of the kernel $F_{j;h,k}$
now leads to 
\begin{subequations}
\begin{equation}
 \label{profile2}
u_j(x,t)=e^{t\Delta}a_j(x)+\frac{P_j(x,t)}{|x|^{d+4}}\, + R_j(x,t),
\end{equation}
where $P_j$ is given by 
\begin{equation}
\label{polynomial}
P_j(x,t)  = \gamma_d \sum_{h,k}
\left( (d+2)x_jx_hx_k - |x|^2 \sigma_{j,h,k}(x) \right) \cc_{h,k}(t),
\end{equation}
with $\sigma_{j,h,k}(x)=\delta_{j,h}x_k+\delta_{j,k}x_h+\delta_{h,k}x_j$.
The remainder $$R_j(x,t)=
-\sum_{h,k} \Bigl(R_{j;h,k}^{(1)}(x,t)+R_{j;h,k}^{(2)}(x,t)\Bigr)$$ is given by~:
\begin{align}
R_{j;h,k}^{(1)}(x,t)&=
|x|^{-(d+1)} \int_0^t {\cal E}_{h,k}(t-s) \:
 \Psi_{j;h,k}\biggl(\frac{x}{\sqrt{s+1}}\biggr)
\,ds,\\
R_{j;h,k}^{(2)}(x,t)&=
\int_0^t v_{h,k}(s)\ast F_{j;h,k}(t-s)\,ds.
\end{align}
\end{subequations}
The functions $\Psi_{j;h,k}$ are given by \eqref{F}.
\begin{remark}
The above remainder is not small when $|x|\leq \sqrt{t}$.
As the solution $u(x,t)$ is smooth
at least for small $t>0$, the homogeneous polynomial and the remainder have to behave
in exactly anti-symmetrical ways when $|x|\to 0$. The same compensation also
occurs for a.e. $x\in\R^d$ when $t\to0$.
\end{remark}

\bigskip
The polynomial profile $\widetilde{P}(x,t)=|x|^{-d-4}P(x,t)$ has no vorticity,
\ie the matrix
$$\operatorname{rot}\widetilde{P} =
\big(\partial_i\widetilde{P}_j-\partial_j\widetilde{P}_i\big)_{i,j}$$
is identically zero. This means that the polynomial profile is
a gradient vector field. In fact, one may check immediately that~:
\begin{equation}
\frac{P(x,t)}{|x|^{d+4}} = \nabla\Pi
\qquad\text{with}\qquad
\Pi(x,t)= \gamma_d\left(\frac{\Tr\cc(t)}{d\,|x|^d}
-\sum_{h,k}\frac{x_h x_k}{|x|^{d+2}}\cdot\cc_{h,k}(t)\right).
\end{equation}

\subsection{General bounds of the remainder terms $R_{j;h,k}$}
\label{sub2.3}
Let us now compute some upper bounds of the remainder terms.
This second part of the proof is also shared by
Theorem~\ref{theorem1} and \ref{theorem3}.

\paragraph{Bound of $R_{j;h,k}^{(1)}(x,t)$.}
The bound~\eqref{bound-Psi} gives~:
$$\Psi_{j;h,k}\biggl(\frac{x}{\sqrt{4(s+1)}}\biggr) \leq C \exp
\biggl(-\frac{c\, |x|^2}{4(s+1)}\biggr),$$
hence
\begin{equation}\label{rem1}
|R_{j;h,k}^{(1)}(x,t)|\leq
 C\, |x|^{-d-1}\exp\biggl(-\frac{c\, |x|^2}{4(t+1)}\biggr)
\:\int_0^t\norme[L^2]{u(s)}^2\,ds.
\end{equation}

\paragraph{Bound of $R_{j;h,k}^{(2)}(x,t)$.}
Since  $\int_{\R^d}v_{h,k}(x,s)\,dx=0$, the second remainder can also be
written~:
$$R_{j;h,k}^{(2)}(x,t)=
\int_0^t\int_{\R^d} v_{h,k}(y,s)\,\left( F_{j;h,k}(x-y,t-s)
-F_{j;h,k}(x,t-s)\right)\,ds\,dy.$$
The Taylor formula gives~:
\begin{align}
|R_{j;h,k}^{(2)}(x,t)| &\leq 
\int_0^t\int_{|y|\leq|x|/2} |y|\,|v_{h,k}(y,s)|\,
\sup_{|z|\leq |x|/2} |\nabla F_{j;h,k}(x+z,t-s)|
\,ds\,dy \notag \\
&\qquad + \int_0^t \left(\int_{|y|\geq |x|/2} |v_{h,k}(y,s)|\,dy\right)
|F_{j;h,k}(x,t-s)|\,ds  \label {R2}\\
&\qquad + \int_0^t\int_{|y|\geq |x|/2} |v_{h,k}(y,s)|\,|F_{j;h,k}(x-y,t-s)|\,ds\,dy. \notag
\end{align}
Thanks to~\eqref{bound-Psi}--\eqref{F}, one has
$|\nabla F_{j;h,k}(x,t)|\le C|x|^{-(d+2)}$
uniformly for $t>0$. Applying~\eqref{L1F} as well, we get~:
\begin{align}
|R_{j;h,k}^{(2)}(x,t)| &\leq 
C \biggl( \int_0^t \!\! \int_{|y|\le |x|/2} |y|\,|v_{h,k}(y,s)|\,dy\,ds\biggr)
|x|^{-(d+2)}\notag \\
&\qquad+\biggl(\int_0^t\!\!\int_{|y|\ge |x|/2} |v_{h,k}(y,s)|\,dy\,ds\biggr)
 |x|^{-(d+1)} \label{rem2}\\
&\qquad+\int_0^t (t-s)^{-1/2}\,\sup_{|y|\ge |x|/2} |v_{h,k}(y,s)|\,ds. \notag
\end{align}

To conclude the proofs of Theorem~\ref{theorem1} and \ref{theorem3},
we shall now use the assumptions on $u$ to estimate \eqref{rem1} and \eqref{rem2}.

\subsection{Local-in-time solutions. Proof of Theorem~\ref{theorem1}}
\label{sub2.4}
The goal of this section is to get upper bounds of the above remainders
that provide valuable information for short time.
In particular, in view of the proof of the lower bounds~\eqref{lower-upper bounds},
it is of interest to have information on the behavior as $t\to0$
of the last term appearing  in the right-hand side of~\eqref{profile}.

\bigskip
The remainder $R^{(1)}_{j;h,k}$ satisfies
\begin{equation}
|R^{(1)}_{j;h,k}(x,t)| \leq \frac{C\, (t+1)^{1/2}}{|x|^{\min\{2\vartheta\,;\,d+2\}}}
\int_0^t \norme[L^2]{u(s)}^2 \,ds.
\end{equation}
Indeed, if $\textstyle{\frac{d+1}{2}<\vartheta\leq\frac{d+2}{2}}$, one has
$$\exp\biggl(-\frac{c\, |x|^2}{4(t+1)}\biggr) \leq C' |x|^{d+1-2\vartheta}
(t+1)^{\vartheta-\frac{d+1}{2}}$$
and if $\vartheta\geq\frac{d+2}{2}$, one also has
$$\exp\biggl(-\frac{c\, |x|^2}{4(t+1)}\biggr) \leq C' |x|^{-1} (t+1)^{1/2}.$$
In both cases, our estimates can blow up
as $t\to\infty$, but not faster than $(1+t)^{1/2}$.

\bigskip
To deal with the remainder $R^{(2)}_{j;h,k}$, one may notice that
the definition of $v$ implies~:
\begin{equation*}
|v_{h,k}(y,s)| \le |u(y,s)|^2+\|u(s)\|_2^2 \,g(y)
\le C(1+|y|)^{-2\ubar\vartheta} \norme[L^\infty_{\ubar\vartheta}]{u(s)}^2
\end{equation*}
with $\ubar\vartheta=\min\{\vartheta\,;\,d+1\}$.
Therefore, since $2\ubar\vartheta>d+1$~:
\begin{align*}
|R_{j;h,k}^{(2)}(x,t)| &\leq 
C \, |x|^{-d-2} \, \int_0^t \norme[L^\infty_{\ubar\vartheta}]{u(s)}^2 \,ds \\
&\qquad+|x|^{-1-2\ubar\vartheta} \left(\int_0^t
 \norme[L^\infty_{\ubar\vartheta}]{u(s)}^2 \,ds \right)\\
&\qquad+(1+|x|)^{-2\ubar\vartheta}\,
\int_0^t\norme[L^\infty_{\ubar\vartheta}]{u(s)}^2 \, (t-s)^{-1/2}\,ds,
\end{align*}
and hence for $|x|\geq1$~:
\begin{equation}
|R_{j;h,k}^{(2)}(x,t)| \leq 
\frac{C\,(t+\sqrt{t})}{|x|^{\min\{2\vartheta\,;\,d+2\}}} \:
\sup_{s\leq t}\norme[L^\infty_{\ubar\vartheta}]{u(s)}^2.
\end{equation}
This ends the proof of~\eqref{profile}.
\medskip
To obtain an asymptotic profile for the pressure,
we need the following simple result on the localization of the derivatives.

\begin{proposition}
\label{space-time deriv}
Given $u\in L^\infty([0,T];L^\infty_\vartheta)$ a solution of
the Navier-Stokes system with Cauchy datum $a=u(0)$ and 
$0\le \vartheta\leq d+1$.
If, for some index~$i$, one has $\partial_i a\in L^\infty_{\vartheta}$,
then~:
\begin{equation}\label{eq:decay_derivatives}
\partial_i u \in L^\infty([0,T];L^\infty_\vartheta).
\end{equation}
If, moreover, $\partial_i a\in L^\infty_\vartheta$ and $\partial_i\partial_j a\in L^\infty_\vartheta$
holds for all $i,j\in\{1,\ldots,d\}$, then 
\begin{equation}
\label{decay 2nd derivatives}
\partial_i\partial_j u \in L^\infty([0,T];L^\infty_\vartheta)
\end{equation}
and
\begin{equation}
\label{decay time derivative}
t^{1/2} \partial_t u \in L^\infty([0,T];L^\infty_\vartheta).
\end{equation}
\end{proposition}

\Proof
Let us first deal with the first order spatial derivatives.
Taking the $i$-th derivative in~\eqref{IE} leads to the affine
fixed point problem~:
$$\partial_i u = \Theta(\partial_i u)$$
with $\Theta=(\Theta_1,\ldots,\Theta_d)$ and
$$\Theta_j w = e^{t\Delta}(\partial_i a_j) 
    - 2\sum_{h,k}\int_0^t F_{j;h,k}(t-s)\ast (u_h w_k)(s) \,ds.$$
Proposition~3 of~\cite{Vig05} implies that 
$\Theta$ is a continuous operator on $X=L^\infty([0,T_0];L^\infty_\vartheta)$, $0\le T_0<T$ and that
$$\norme[X]{\Theta(w-w')} \leq
C_0\,T_0^{\frac{1}{2}}\,\sup_{t\in[0,T]}\norme[L^\infty_\vartheta]{u(t)}
\norme[X]{w-w'}.$$
One may therefore choose $T_0>0$ such that $\Theta$ is a contraction of the
Banach space $X$.  Its only fixed point $w=\partial_i u$ belongs
therefore to this function space.
The same argument also holds on $[T_0,2T_0]$,\ldots{} and leads finally
to $\partial_i u \in L^\infty([0,T];L^\infty_\vartheta)$.

Conclusion~\eqref{decay 2nd derivatives} also follows from the contraction mapping theorem
in a similar way.

\medskip
For the time derivative, the starting point is again an
identity that directly follows from~\eqref{IE}, namely,
$$\partial_t u =  \tilde\Theta(\partial_t u),$$
with
$$\tilde\Theta_j(w) = \Delta a_j- \sum_{h,k} F_{j;h,k}(t)\ast (a_ha_k)-
2\sum_{h,k}\int_0^t F_{j;h,k}(s) \ast (u_h w_k)(t-s) \,ds.$$
The Banach space we deal with is
$$Y=\{w\,;\,t^{1/2}\norme[L^\infty_\vartheta]{w(t)} \in
L^\infty([0,T_0])\}.$$
Proposition~3 of~\cite{Vig05} implies now that
$$-\Delta a_j+ \sum_{h,k} F_{j;h,k}(t)\ast (a_ha_k)\in Y$$
and
$$\big\|\tilde\Theta(w-w')\big\|_{Y} \leq
\pi C_0\,T_0^{\frac{1}{2}}\,\sup_{t\in[0,T]}\norme[L^\infty_\vartheta]{u(t)}
\norme[Y]{w-w'}.$$
Here, we have used the fact that~:
$$\forall t>0, \qquad \int_0^{t}\frac{ds}{\sqrt{s(t-s)}}=\pi.$$
The conclusion now follows in the same lines as above.
\endProof

We can now establish~\eqref{profile pressure}~:
The pressure is defined up to an arbitrary
function of $t$ by~:
$$-\nabla p = (\partial_t-\Delta)u+\div(u\otimes u).$$
Let us now replace $u$ by 
the its profile given by~\eqref{profile2}, that is
$u=e^{t\Delta}a+\nabla\Pi+R$.
One gets~:
$$-\nabla p = \nabla(\partial_t-\Delta)\Pi+(\partial_t-\Delta)R+\div(u\otimes u).$$
This yields~:
\begin{equation}
\label{profile pressure2}
p(x,t)=p_0-\gamma_d
\left(\frac{\Tr{\cal E}(t)}{d\,|x|^d}
-\sum_{h,k}\frac{x_h x_k}{|x|^{d+2}}\cdot{\cal E}_{h,k}(t)\right) +q(x,t),
\end{equation}
where the remainder term $q(x,t)$ satisfies
\begin{equation}
-\nabla q = -\nabla\Delta\Pi + (\partial_t-\Delta)R+\div(u\otimes u).
\end{equation}

\noindent%
Let us show that, for all $t>0$, we have 
$\nabla q = {\cal O}_t(|x|^{-\min\{2\vartheta,d+2\}})$.
\begin{itemize}
\item[--]
We obviously have 
$\nabla\Delta\Pi=
{\cal O}_t(|x|^{-\min\{2\vartheta,d+2\}})$,
since the left hand side is  a homogeneous function
of degree  $-d-3$ which is smooth for $x\not=0$.
\item[--]
The term $\div(u\otimes u)=(u\cdot\nabla)u$ belongs uniformly
to $L^\infty_{2\vartheta}$ because
of~\eqref{eq:decay_derivatives}.
\item[--]
The remainder is the sum of two terms: one checks immediately
that $(\partial_t-\Delta)R^{(1)}$ is exponentially decaying as $|x|\to\infty$.
The second term is
$$(\partial_t-\Delta)R_{j;h,k}^{(2)}(x,t) = 
v_{h,k}(0)\ast F_{j;h,k}(t)+
\int_0^t (\partial_t-\Delta)v_{h,k}(t-s)\ast F_{j;h,k}(s)\,ds$$
where $v_{h,k}(x,t)=u_hu_k - {\cal E}_{h,k}(t)\,g(x)$ and
$t^{1/2}\partial_t v_{h,k}$ belongs to $L^\infty([0,T];L^\infty_{2\vartheta})$.

We now use again that $\displaystyle\int v_{h,k}(0)\,dx=0$ :
if we apply Lemma~\ref{lemma classical}, the computation~\eqref{rem2}
shows that the first term is bounded in 
$t^{-1/2} L^\infty([0,T];L^\infty_{\min\{2\vartheta;d+2\}})$.
On the other hand, $\partial_t v_{h,k}$ and $\Delta v_{h,k}$ also have
a vanishing integral.
Therefore, using the estimates on the space-time derivatives provided by
Proposition~\ref{space-time deriv} shows that the second 
term belongs to $L^\infty([0,T];L^\infty_{\min\{2\vartheta;d+2\}})$.
\end{itemize}
Hence we get $\nabla q\in {\cal O}_t(|x|^{-\min\{2\vartheta,d+2\}})$.
Our last step is the following elementary estimate~: 

\begin{lemma}
\label{elementary lemma}
Let $\alpha>1$ and $f\in C^1(\R^d)$ such that $\nabla f\in L^\infty_\alpha$.
Then there is a constant $c$ such that $f-c\in L^\infty_{\alpha-1}$. 
\end{lemma}

\Proof
For any $\omega\in \R^d$, $|\omega|=1$, let
$\ell_\omega\equiv\lim_{r\to\infty} f(r\omega)
  =f(0)+\displaystyle\int_0^\infty \nabla f(s\omega)\cdot\omega\,ds$.
If $\tilde\omega$ is another point of the unit sphere then for all $r>0$ we have
$$|\ell_\omega-\ell_{\tilde\omega}|\le \int_r^\infty |\nabla f(s\omega)|\,ds
   +C r\sup_{|x|\ge r}|\nabla f(x)|
   + \int_r^\infty |\nabla f(s\tilde\omega)|\,ds.$$
Letting $r\to\infty$ we get that $c\equiv \ell_\omega$
is independent of $\omega$.
But
$$|f(r\omega)-c|\le \int_r^\infty |\nabla f(s\omega)|\,ds\le C(1+r)^{-\alpha+1}$$
and the conclusion follows.
\endProof

\medskip
The standard properties of strong solutions imply that $q(x,t)$ is smooth for $x\not=0$.
Applying this lemma (for fixed~$t$, $0<t<T$), with $f(x)=\chi(x)q(x,t)$, where $\chi$
is a smooth function such that $\chi(x)\equiv0$ for $|x|\le r$ and $\chi\equiv1$ for $|x|\ge r'$
for some $0<r<r'$ implies that,
$q(x,t)= c+{\cal O}_t(|x|^{-\min\{2\vartheta-1,d+1\}})$.
This completes the proof of Theorem~\ref{theorem1}.
\endProof

\bigskip
For later use, let us note explicitly that if $R(x,t)={\cal O}_t(|x|^{-\min\{2\vartheta,d+2\}})$
denotes the last term in the right hand side of~\eqref{profile}, then we
have proved that, for all $0\le t\le T$,
\begin{equation}
\label{bound on R}
|R(x,t)|\le \frac{C(\sqrt t+t)\,\|a\|_{L^\infty_\vartheta}}{|x|^{\min\{2\vartheta,d+2\}} } .
\end{equation}

\bigskip
\subsection{Global-in-time solutions. Proof of Theorem~\ref{theorem3}}
\label{sub2.5}
Let us now focus on long time asymptotics.
%The upper bound~\eqref{rem1} suggests that $R^{(1)}_{j;h,k}(x,t)$ is exponentially
%decaying when $|x|\to+\infty$, but that this term does not get smaller when $t\to+\infty$.
%It should therefore be included in the asymptotic profile of $u$.
Let $u$ be a global solution satisfying~\eqref{assumption_largetime}.
Going back to~\eqref{profile2}, we see that the  profile~\eqref{detailed profile}
holds with
\begin{equation*}
E_j(y,t) = - \sum_{h,k} \int_0^t\!\!\int_{\R^d}
(u_hu_k)(y,t-s)\,\Psi_{j;h,k}\!\left(\sqrt{\frac{t+1}{s+1}}\:y\right)
\,dz\,ds
\end{equation*}
and
$${\cal R}_j(x,t)=-\sum_{h,k} R^{(2)}_{j;h,k}.$$

%
%\begin{gather}
%\label{bound_Ej}
%|E_j(y,t)|\leq C' e^{-c|y|^2} \norme[{L^2([0,t]\times\R^d)}]{u}^2 \\[2pt]
%\intertext{and, for $\sigma=\vartheta-(d+2)/2$, $|y|\geq 1$ and $\epsilon>0$~:}
%\label{bound_epsilon}
%|\varepsilon_j(y,t)|\leq C_\epsilon \, |y|^{-(2\sigma-\epsilon)} \, (1+t)^{-\sigma}.
%\end{gather}
%\end{theorem}
%

The bound~\eqref{bound-Psi} immediately implies~\eqref{double bound}.
To prove~\eqref{remainder_asympt}, we start by observing that
assumption~\eqref{assumption_largetime} implies~:
$$\norme[L^2]{u(t)}^2 \leq C_\epsilon \,
(1+t)^{-\left(\vartheta-\frac{d}{2}-\epsilon\right)}$$
for any $\epsilon>0$ and therefore, letting $v=(v_{h,k})$,
$$ |v(y,s)| \le |u(y,s)|^2+\|u(s)\|_2^2 \,g(y)
\leq C \, (1+|y|)^{-2\alpha}  \, (1+s)^{-(\vartheta-\alpha)}$$
for $d/2<\alpha\leq\vartheta$.

A consequence of~\eqref{F} is that, for all $0\le \beta\le d+1$ and $0\le \gamma\le 1$
\begin{align*}
&|F_{j;h,k}(x,t)|\le C|x|^{-\beta}t^{-(d+1-\beta)/2},\\
&|\nabla F_{j;h,k}(x,t)|\le |x|^{-(d+1+\gamma)}t^{-(1-\gamma)/2}.
\end{align*}

Therefore, coming back to \eqref{R2} we get
\begin{align}
|{\cal R}(x,t)| &\leq 
C \biggl(\int_0^t \!\! \int_{|y|\le |x|/2}
(1+|y|)^{1-2\alpha}(1+s)^{-\vartheta+\alpha}(t-s)^{-\frac{1}{2}+\frac{\gamma}{2}}\,dy\,ds\biggr)
|x|^{-(d+1+\gamma)}\notag \\
&\qquad+C\biggl(\int_0^t\!\!\int_{|y|\ge |x|/2}    (1+|y|)^{-2\alpha}(1+s)^{-\vartheta+\alpha}|x|^{-\beta}
    (t-s)^{-\frac{d+1-\beta}{2}}  \,dy\,ds\biggr) \\
&\qquad+C\int_0^t (t-s)^{-1/2} (1+|x|)^{-2\alpha}(1+s)^{-\vartheta+\alpha}\,ds, \notag
\end{align}
where, in the last integral, we have also used Lemma~\ref{lemma classical}.

Let us call $I_1$, $I_2$ and $I_3$ the three terms of the right-hand side.
To estimate  $I_1$, we fix a small $\epsilon>0$ and choose $\alpha=\frac{d+1}{2}+\epsilon$.
Then we write $I_1=I_{1,1}+I_{1,2}$, where these two terms are obtained by splitting
the integral $\int_0^t$ into $\int_0^{t/2}$ and $\int_{t/2}^t$.
Then we have, for all $t\ge1$,
$$ I_{1,1}\le C_\gamma\,|x|^{-(d+1+\gamma)}\cdot
\begin{cases}
t^{-\frac{1}{2}+\frac{\gamma}{2}}, &\hbox{if $\vartheta>\frac{d+3}{2}$}\\
t^{\frac{d+2+\gamma}{2}-\vartheta} &\hbox{if $\frac{d+1}{2}<\vartheta\le \frac{d+3}{2}$}
\end{cases}
$$
and
$$I_{1,2}\le C_{\gamma,\epsilon}\, |x|^{-(d+1+\gamma)}t^{\frac{d+2+\gamma}{2}+\epsilon-\vartheta}.$$
Thus,
\begin{equation}
\label{bound I1}
I_1\le  C_{\gamma,\epsilon}\,|x|^{-(d+1+\gamma)}\cdot
\begin{cases}
 t^{-\frac{1}{2}+\frac{\alpha}{2}}, &\hbox{if $\vartheta>\frac{d+3}{2}$}, \\
 t^{\frac{d+2+\alpha}{2}-\vartheta+\epsilon}, &\hbox{if $\frac{d+1}{2}<\vartheta<\frac{d+3}{2}$},
\end{cases}
\end{equation}

To estimate $I_2$, we choose again $\alpha=\frac{d+1}{2}+\epsilon$ and $\beta=d+\gamma$.
Then the same argument as before shows that $I_2$ can be bounded as in~\eqref{bound I1}.
To estimate $I_3$, we take  $\alpha=\frac{d+1+\gamma}{2}$. This choice shows 
that $I_3$ is also bounded by the function on the right hand side of~\eqref{bound  I1}.
Summing all these bounds completes the proof of Theorem~\ref{theorem3}.

\endProof

\subsection{Criterion for the vanishing of $\nabla\Pi$ -- Proof
of Proposition~\ref{theorem2}}
\label{sub2.6}

We now prove Proposition~\ref{theorem2}, which we restate in a more
complete form.

\begin{proposition}
\label{proposition-CNS}
For any real matrix $\cc=(\cc_{h,k})$, let us define a family
of homogeneous polynomials by
\begin{equation}\label{polyQ}
Q_j(x) = \sum_{h,k}
\left(|x|^2 \sigma_{j,h,k}(x) - (d+2)x_jx_hx_k\right) \cc_{h,k}.
\end{equation}
The following assertions are equivalent~:
\begin{enumerate}
\item
The matrix is proportional to the identity matrix, \ie
\begin{equation}
\label{orth2}
\forall h,k\in\{1,\ldots,d\},\qquad
\cc_{h,k}=\alpha\,\delta_{h,k}
\end{equation}
with $\alpha=\frac{1}{d}\Tr\cc$.
\item
$Q_j\equiv 0$ for all indices $j\in\{1,\ldots,d\}$.
\item
There exists an index $j\in\{1,\ldots,d\}$ such that $Q_j\equiv 0$.
\item
There exists an index $j\in\{1,\ldots,d\}$ such that $\partial_j Q_j\equiv 0$.
\end{enumerate}
\end{proposition}

\medskip
Putting the terms $x_jx_\ell^2$ in factor in \eqref{polyQ}, one gets the
following expression for the $j^\text{th}$ component of $Q$~:
\begin{equation*}\begin{split}
Q_j(x)= x_j\sum_{\ell=1}^d \left\{\Tr\cc-d\cc_{\ell,\ell}
+2(\cc_{j,j}-\cc_{\ell,\ell})\right\} x_\ell^2 &\\
+ 2|x|^2 \tilde{\cc}(\e{j},x)-(d+2)x_j & \tilde\cc(x,x)
\end{split}\end{equation*}
where $\e{i}=(0,\ldots,0,1,0,\ldots,0)$ denotes the canonical basis of $\R^d$
and $\tilde\cc$ is the bilinear form defined by the non-diagonal
coefficients of $\cc$~:
$$\tilde\cc(u,v)=\sum_{h\neq k} \cc_{h,k} u_h v_k.$$

Relations \eqref{orth2} express the fact that the 
matrix $\cc=(\cc_{h,k})_{1\leq h,k \leq d}$ is a scalar multiple
of the identity matrix.
In such a case, one can immediately check on the previous
expression that~$Q_j(x)=0$.

Let us prove conversly that $\partial_jQ_j\equiv0$
implies~$\cc=\alpha\id$. One has~:
$$\begin{aligned}\partial_jQ_j(x)&=
\sum_{\ell=1}^d \left\{(1+2\delta_{j,l})(\Tr\cc-d\cc_{\ell,\ell})
+2(\cc_{j,j}-\cc_{\ell,\ell})\right\} x_\ell^2\\
&\qquad - 2 d x_j \tilde{\cc}(\e{j},x)-(d+2)\tilde\cc(x,x).
\end{aligned}$$
The fact that $\partial_jQ_j(\e{i})=0$ for all $i$ implies
$$\forall \ell\in\{1,\ldots,d\},\qquad
(1+2\delta_{j,l})(\Tr\cc-d\cc_{\ell,\ell})
+2(\cc_{j,j}-\cc_{\ell,\ell})=0$$
and hence $\cc_{i,i}=\frac{1}{d}\Tr\cc$ ($i=1,\ldots,d$),
\ie all the diagonal entries of $\cc$ are equal.
Therefore~:
$$\partial_jQ_j(x)=- 2 d x_j
\tilde{\cc}(\e{j},x)-(d+2)\tilde\cc(x,x)$$
and this expression should vanish identically.
A new derivation with respect to $x_j$ gives
$$\partial_j^2Q_j=-4(d-1)\tilde\cc(\e{j},x)=0,$$ \ie $\tilde\cc(\e{j},x)=0$
as~$d\geq2$, and hence $\tilde\cc(x,x)\equiv0$. This proves that
the matrix $\cc$ is a scalar multiple of the identity matrix.

\endProof

\section{Applications}
\label{section3}

%Before going on to the proof of the above results,
%let us explore several of their consequences.
Let us now explore a few consequences of the above results.

\subsection{Instantaneous spreading property}
\label{sub3.1}

It is a consequence of the result of~\cite{BraM02} that, if the components of the
initial data have no special symmetries, then
the corresponding solution $u(x,t)$ of \NS satisfies 
\begin{equation} \label{BraMres}
\liminf_{R\to\infty}R \int_{R\le|x|\le 2R} |u(x,t)|\,dx >0\,,
\end{equation}
for $t>0$ belonging at least to a sequence of points $t_k$ converging to zero as $k\to\infty$.
In particular for those~$t$,
one has~:
$$\int_{\R^d}|x|\,|u(x,t)|\,dx =+\infty
\qquad\text{and}\qquad\int_{\R^d}|x|^{d+2}|u(x,t)|^2\,dx=+\infty.$$
(See also, \eg\cite{Lem02}, Theorem~25.2). 
The precise assumption  guaranteeing \eqref{BraMres} is  the non-orthogonality
of the components with respect to the $L^2$-inner
product, \ie  one can find $j\neq k$ in $\{1,\ldots,d\}$
such that
\begin{subequations}\label{non-orth}
\begin{equation}\label{non-orthA}
\int_{\R^d}a_j(x)\,a_k(x)\,dx\neq0
\end{equation}
or such that
\begin{equation}\label{non-orthB}
\int_{\R^d}a_j^2(x)\,dx\neq\int_{\R^d}a_k^2(x)\,dx.
\end{equation}
\end{subequations}
Even if \eqref{BraMres} already explains that the
limitation $\ubar\vartheta\leq d+1$
in Theorem~\ref{theorem-intro} is optimal for generic flows,
such a condition does not provide much information on the pointwise decay of $u$,
as $|x|\to\infty$.

\bigskip
Theorem~\ref{theorem1} and Proposition~\ref{theorem2} not only provide
such information, but also allow us to give a simpler proof of these facts.
The instantaneous spreading property is fully described by
the following result.

\begin{theorem}\label{corollary1}
For $\vartheta>d+1$, let $a\in L^\infty_\vartheta$ be a divergence-free vector field.
Let $u$ be the corresponding solution of \NS
in ${\cal C}_w\left([0,T];L^\infty_{d+1}\right)$.
For $0<t\le T$, we set 
$$\kappa_t=\max\{1,t^{-1/2},t^{-1/(\vartheta-d-1)}\}.$$
\begin{enumerate}
\item
There is a constant $c>0$ such that, for $0<t\le T$ and $|x|\ge c\kappa_t$~:
\begin{equation}\label{upper}
|u(x,t)| \leq c \, t \, |x|^{-(d+1)}.
\end{equation}
\item
Conversely, if \eqref{non-orth} holds for a couple of indices $(j,k)$,
then there exists $t_0\in(0,T]$ and a constant $c'>0$ such that
for all $0<t\le t_0$ and all $x$ in a conic neighborhood of the~$x_j$ or $x_k$ axis,
with $|x|\ge c\kappa_t$~:
\begin{equation}\label{lower-j}
|u_j(x,t)|\ge c'\,t\,  |x|^{-(d+1)}.
\end{equation}

\item
Actually, if \eqref{non-orth} holds, the lower bound~\eqref{lower-j} holds in almost 
all directions~: the set
\begin{equation}\label{lower-set}
\Sigma = \left\{ \sigma\in\mathbb{S}^{d-1} \,;\,
\underset{{\substack{t\to 0^+\\|x|\to\infty \\ x\in\R\sigma}}}{\operatorname{lim\, inf}}
\left( t^{-1} \, |x|^{d+1} \, |u_j(x,t)| \right) = 0
\right\}
\end{equation}
is a closed subset of the sphere $\mathbb{S}^{d-1}$,  of measure zero.
\end{enumerate}
\end{theorem}

\begin{remark}
The non-orthogonality assumption~\eqref{non-orth} cannot be removed.
Indeed, \cite{BraM02} contains an example
of symmetric and highly localized flow for which 
\eqref{BraMres} and \eqref{lower-j} break down.
For a further understanding of how symmetries usually lead to highly
localized flows, one may also refer to \cite{HeM03} or \cite{Bra04ii}.
\end{remark}

The previous results show that, for generic flows,
\emph{each component} of the velocity field
instantaneously spreads out in \emph{any direction}.
In particular, one  has the following corollary.
\begin{corollary}
Under the same assumptions as in Theorem~\ref{corollary1}
there exist a time $t_0>0$ and two positive constants $c_0$, $C$ such that for all $t\in(0,t_0]$~:
\begin{equation}\label{lower-integral}
 \int_{R\le|x|\le 2R} |u_j(x,t)|\,dx \ge \frac{c_0t}{R}
\end{equation}
for all $R\ge C\, t^{-1/\min\{2\,;\,\vartheta-d-1\}}$.
In particular, for $j,k=1,\ldots,d$,
\begin{equation} \label{lower-integral2}
 \int_{\R^d}|x_k|^{\vartheta p}|u_j(x,t)|^p\,dx=+\infty,\qquad
\end{equation}
as soon as $1\le p<+\infty$ and $\vartheta+\frac{d}{p}\ge d+1$.
\end{corollary}

\begin{remark}
The spreading property~\eqref{lower-integral2} should be compared with
F. Vigneron's result \cite{Vig05}, according to which
$$\int_{\R^d}(1+|x|)^{\vartheta p}|u(x,t)|^p\,dx \le C_t \int_{\R^d}
(1+|x|)^{\vartheta p}|a(x)|^p\,dx$$
as long as the solution exists in $L^p(\R^d)$, provided that
$p>d$ and $\vartheta+\frac{d}{p}<d+1$.
\end{remark}

\Proof
First, let us establish the upper bound~\eqref{upper}.
Theorem~\ref{theorem1}, together with~\eqref{bound on R}, immediately implies, for $|x|\ge1$~:
$$|u(x,t)-e^{t\Delta}a(x)|\le C\,t\,|x|^{-(d+1)} + C\sqrt t\, |x|^{-(d+2)}\le C\,t\,|x|^{-(d+1)},$$
for all $x\in\R^d$ such that $|x|\ge t^{-1/2}$.
Moreover, if $|x|\geq (2t^{-1})^{1/(\vartheta-d-1)}$, one has~:
\begin{equation}
\label{bound heat}
|e^{t\Delta}a(x)|\leq C |x|^{-\vartheta} \leq \frac{C}{2} t |x|^{-d-1}.
\end{equation}
This proves~\eqref{upper}.
Let us now focus on the lower bound \eqref{lower-j} of $|u_j(x,t)|$.

\medskip
Let $j\not=k$ such that $\alpha\equiv\int_{\R^d}(a_ja_k)\,dx\not=0$.
Then, for some $t_0>0$, possibly  depending on~$a$, and all $0<t\le t_0$,
we have $$|\cc_{j,k}(t)|\ge \frac{|\alpha| t}{2}.$$
Let $\epsilon>0$ and let $\Gamma_k=\{x \,;\, |x_r|< \epsilon |x_k| \enspace(r\neq k)\}$
be a conical neighborhood of the $x_k$-axis.
Recalling~\eqref{polynomial} we get, for $\epsilon$ small and $R$ large enough~:
$$\forall x\in\Gamma_k, \qquad |x|\geq R \quad\Longrightarrow\quad
|P_j(x,t)|\ge \frac{|\alpha| t}{3}|x_k|^3.$$
Using~\eqref{profile2}, \eqref{bound on R} and the first of~\eqref{bound heat} now leads,
for large enough $|x|$ and $x\in\Gamma_k$, to
$$|u(x,t)|\ge |u_j(x,t)|\ge \frac{|\alpha|t}{4} |x_k|^{-(d+1)}\ge 
 \frac{|\alpha|t}{4} |x|^{-(d+1)}$$ and  \eqref{lower-j} follows
in this case.

In the second case we choose $j\not=k$
such that $\int_{\R^d}a_j^2\,dx\not=\int_{\R^d}a_k^2\,dx$.
Then we have $$d\int_{\R^d} a_j^2\,dx\not=\sum_{m=1}^d \int_{\R^2} a_m^2\,dx$$
(otherwise, $d\int a_k^2\,dx\not=\sum_{m=1}^d \int_{\R^d}a_m^2\,dx$, and
we should exchange $j$ with $k$). 
Let us set $$\beta=\sum_{m=1}^d \int_{\R^d}a_m^2\,dx-d\int_{\R^d}a_j^2\,dx.$$
Arguing as before, we see that there exists a conic neighborhood $\Gamma_j$
of the $x_j$ axis,  such that for all $x\in \Gamma_j$ and $|x|$ large enough
we have $$|u(x,t)|\ge |u_j(x,t)|\ge \frac{|\beta|t}{4} |x_j|^{-(d+1)}\ge 
 \frac{|\beta|t}{4} |x|^{-(d+1)}.$$
Then  \eqref{lower-j} follows in this second case as well.

\smallskip
Let us now prove the last statement of Theorem~\ref{corollary1}.
The map $\displaystyle s\mapsto\int_{\R^d}(u_hu_k)(x,s)\,dx$ is continuous.
Therefore,~\eqref{polynomial} implies that
$$\forall x\in\R^d, \qquad \lim_{t\to0}\frac{1}{t} P(x,t) = {\cal P}(x)$$
where ${\cal P}=({\cal P}_1,\ldots{\cal P}_d)$ is given by
$${\cal P}_j(x) = \gamma_d \, \sum_{h,k}
\left(\int_{\R^d}a_ha_k\right)
\left((d+2)x_jx_hx_k  - |x|^2 \sigma_{j,h,k}(x)  \right),$$
which is a homogeneous polynomial of degree exactly three.
According to Proposition~\ref{theorem2}, the assumption~\eqref{non-orth}
means that ${\cal P}_j\not\equiv0$ for all $j=1,\ldots,d$.

The convergence of $t^{-1}\,P(x,t)$ to ${\cal P}(x)$ is uniform
when $x$ belongs to the (compact) unit sphere.  Let us define
the following dense open subsets of $\Sphere^{d-1}$~:
$$\Omega_j=\{\omega \in \Sphere^{d-1}\colon {\cal P}_j(\omega)\not=0\}.$$
% For example, when $d=2$, $\Omega_j$ consists of the whole $\Sphere^1$ but a
% set of at most six points (this is seen e.g. writing ${\cal P}(x)$
% in polar coordinates).
Given $\omega\in \Omega_j$, let us define $T_{\omega}>0$ as the supremum of
$t\leq T'$ such that $$\frac{1}{t}|P_j(\omega,t)|\geq
\frac{1}{2} |{\cal P}_j(\omega)|.$$
Also let $c_{\omega}=\frac{1}{4} |{\cal P}_j(\omega)|$.

\smallskip
From~\eqref{profile2},~\eqref{bound on R} and the obvious estimate
 $|e^{t\Delta}a(x)|\le C|x|^{-\vartheta}$,
we get,
 for $\omega=x/|x|\in\Omega_j$~:
$$ |u_j(x,t)|\ge 2c_{\omega}\,t\,|x|^{-(d+1)}- C\,t^{1/2}|x|^{-(d+2)}-C|x|^{-\vartheta}\ge c_{\omega}\,t\, |x|^{-(d+1)}$$
for all $0\le t\le T_{\omega}$ and $|x|\ge C_{\omega}\,t^{-\vartheta^*}$.
The complement of $\Omega_j$ is an algebraic surface
and therefore the set $\Sigma=\mathbb{S}^{d-1}\backslash\Omega_j$
has measure zero in $\mathbb{S}^{d-1}$.

\bigskip
Finally, \eqref{lower-integral} and the corollary follow
immediately from~\eqref{lower-j} and the fact that a function bounded
from below, at infinity, by $|x|^{-d-1}$
does not belong to any weighted Lebesgue space $L^p(\R^d,(1+|x|)^{\vartheta p}\,dx)$ when $1\leq p<\infty$
and $\vartheta+\frac{d}{p}\geq d+1$.

\endProof

\begin{remark}
The set of exceptional directions where, at some time $t>0$,
our lower bound~\eqref{lower-j}
breaks down corresponds to the roots of the polynomial $P(\cdot,t)$
on the sphere $\Sphere^{d-1}$.
In dimension two we have a more precise description
of such a set~: it is made of at most {\it six directions\/}.
This is due to the fact that a homogeneous polynomial in $\R^2$
of degree exactly three has at most six zeros on $\Sphere^1$,
as is easily checked by passing to polar coordinates.
In principle, this conclusion is valid only during a short time interval $[0,t_0]$.
Indeed we cannot exclude that after some time there are flows featuring
some kind of creation of symmetry: it may happen that, for some time $t_1>t_0,$
$\int u_1(x,t_1)^2\,dx=\int u_2(x,t_1)^2\,dx$ and $\int (u_1u_2)(x,t_1)\,dx=0$,
in a such way that the velocity field first instantaneously spreads out but then
recovers a good localization at time $t_1$. 
Subsequently, the flow would remain localized, 
or spread out again, depending on whether the components $u_1$
and $u_2$ remain orthogonal or not, after $t_1$.
However, no example of such a somewhat pathological flow is known so far.
\end{remark}

\subsection{Lower bounds of solutions in weighted spaces}
\label{sub3.2}

Let us establish a few consequences of Theorem~\ref{theorem3}.
Throughout this section we suppose $t\ge1$.
For well localized data, {\eg}  when $a\in L^\infty_\vartheta$ 
with $\vartheta>\frac{d+2}{2}$, the limit
$$\nabla\Pi_\infty(x)\equiv\lim_{t\to\infty}\nabla\Pi(x,t)$$
is well defined.
In this case, a consequence of~\eqref{detailed profile} is that, for some $\beta>0$,
\begin{equation}
\label{large time}
\Bigl|u(x,t)-e^{t\Delta}a(x)-\nabla\Pi_\infty(x)\Bigr|
\le C|x|^{-d-1}e^{-|x|^2/(t+1)}+C|x|^{-d-1}\,t^{-\beta}.
\end{equation}
This leads us to introduce, for all $A>0$, the region
$$\mathcal{D}_A(t)=\{x\in\R^d \,;\, |x|^2\geq A(t+1)\}.$$
Since, generically, $\nabla\Pi_\infty\not\equiv0$,
several  lower bounds for the large time behavior
of $u$ can be obtained as an easy consequence of~\eqref{large time}.
For example, if we introduce the weighted norm
\begin{equation*}
\norme[{L^p_\alpha}]{f}=\biggl(\int |f(x)|^p(1+|x|)^{p\alpha}\,dx\biggr)^{1/p},
\end{equation*}
then, taking $A>0$ large enough,  for all $1\le p<\infty$ and $\alpha\ge0$, such that
\begin{equation}
\label{restriction}
\alpha+\frac{d}{p}<d+1,
\end{equation}
we get, for all $t>0$ large enough,
\begin{equation}
\label{large time bounds}
\begin{aligned}
\norme[{L^p_\alpha}]{u(t)-e^{t\Delta}a}^p 
     &\ge \int_{\mathcal{D}_A(t)} |u(x,t)-e^{t\Delta}a(x)|^p(1+|x|)^{p\alpha}\,dx\\
    &\ge \frac{1}{2}\int_{\mathcal{D}_A(t)} |\nabla \Pi_\infty(x)|^p(1+|x|)^{p\alpha} \,dx\\
    &\ge C\,(At)^{-\frac{p}{2}(d+1-\alpha-\frac{d}{p})}.
\end{aligned}
\end{equation}

If the datum is highly oscillating (for example, if the  Fourier transform of $a$ satisfies some suitable
vanishing condition at the origin) then the $L^p_\vartheta$ norms of $e^{t\Delta}a$ decay faster as $t\to\infty$
than the right hand side of~\eqref{large time bounds}.
Then, \eqref{large time bounds} will be in  fact a lower bound for~ $\norme[{L^p_\alpha}]{u(t)}^p $
in this case.
A similar conclusion remains true if we drop this assumption on the oscillations
and start with a datum that is simply well localized.
Indeed, we have the following~:

\begin{corollary}
\label{corollary4}
Let $u$ be as in Theorem~\ref{theorem3}, starting from $a\in L^\infty_\vartheta$, with
$\vartheta>d+1$.
We also assume that $\nabla\Pi_\infty\not\equiv0$.
Then there exist $t_0>0$ and  a constant $c>0$ such that, for all  $1\le p<\infty$ and $\alpha\ge0$, satisfying~\eqref{restriction},
we have, for all $t\ge t_0$~:
\begin{equation}
\label{large time lower2}
\norme[{L^p_\alpha}]{u(t)} \geq
c\,t^{{-\frac{1}{2}(d+1-\alpha-\frac{d}{p})}}. 
\end{equation}
Moreover, for all  $0\le\alpha\le d+1$ and all~$t\ge t_0$~:
\begin{equation}
\label{large time lower3}
\norme[{L^\infty_\alpha}]{u(t)} \ge c\,t^{-\frac{1}{2}(d+1-\alpha)}.
\end{equation}
\end{corollary}

The lower bound~\eqref{large time lower2}
 was already known for $p=2$ and $0\le \alpha\le 2$
(see, \eg \cite{Sch91}, \cite{BaJ05}), or~$1\le p\le \infty$ and $\alpha=0$
(see \cite{FujM01}).
The decay profile~\eqref{assumption_largetime}, under the assumption
of Corollary~\ref{corollary4}, immediately implies the (slightly weaker) upper bound
$\norme[{L^p_\alpha}]{u(t)} \leq
c_\epsilon\,t^{{-\frac{1}{2}(d+1-\alpha-\frac{d}{p}-\epsilon)}}$
for all~$\epsilon>0$.
In fact, the ``sharp'' upper bound (\ie the bound with $\epsilon=0$)
has been obtained, at least for $p\ge2$ and with some
additional restrictions on $\alpha$, by many authors
(see \cite{KuT05} and the references therein).

\medskip
\Proof
By our assumptions,  $a\in L^1(\R^d)$ and $\hbox{div\,}a=0$. Thus,
$\displaystyle\int a(y)\,dy=0$.
A direct computation (using the same method as in the proof of~\eqref{R2})
then yields
$$ |e^{t\Delta}a(x)| \le C(1+|x|)^{-\vartheta}\,(1+t)^{(\vartheta-d-1)/2}.$$
It then follows that, for $t\ge 1$,
$$ \int_{\mathcal D_A(t)} |e^{t\Delta}a(x)|^p(1+|x|)^{\alpha p}\,dx
    \le C\,A^{-\frac{p}{2}(\vartheta-\alpha-\frac{d}{p})}\,t^{-\frac{p}{2}(d+1-\alpha-\frac{d}{p})}.$$
Here the exponent of~$A$ is strictly smaller than that of~\eqref{large time bounds}.
If $A$ is large enough, then a comparison between this inequality and~\eqref{large time bounds}
gives~\eqref{large time lower2}.
The proof of \eqref{large time lower3} is essentially the same.
\endProof

\subsection{Flows with anisotropic decay in the whole space}
\label{section-anisotropic}

This short section contains a positive and a negative result about
flows in $\R^d$ with anisotropic decay at infinity.

\bigskip
Theorem~\ref{theorem1} implies that \NS flows may inherit the anisotropic
decay properties of the initial data, as long as these properties do not violate
the instantaneous spreading limit given by Theorem~\ref{corollary1}.

\begin{proposition}\label{anisotropic-stability}
Let $a$ be a bounded divergence-free vector field and  $u$ the corresponding
solution of \NS in ${\cal C}_w\left([0,T);L^\infty\right)$.  
Let us also assume that there exists a function $m$ such that
\begin{equation}\label{assume_anisotrope}
|e^{t\Delta}a(x)| \leq C_t (1+|x|)^{-\vartheta} \, m(x)^{-1}
\end{equation}
with $\frac{d+1}{2}<\vartheta\leq d+1$ and $1\leq m(x)\leq C
(1+|x|)^{d+1-\vartheta}$.
Then, for all $T'<T$ there exists a constant $C_{T'}$ (this also might depend
on the data) such that~:
\begin{equation}
|u(x,t)| \leq C_{T'} (1+|x|)^{-\vartheta} \, m(x)^{-1}
\end{equation}
for all $t\in[0;T']$.
\end{proposition}
\Proof
This is an obvious consequence of~\eqref{profile} and~\eqref{bound on R}.

\endProof

Let us give some examples of anisotropic weights satisfying~\eqref{assume_anisotrope}.
A \emph{Peetre-type weight} is a measurable function $m:\R^d\to[1;+\infty)$
such that
\begin{equation}\label{peetre}
\exists C_0>0,\quad \forall x,y\in\R^d,\qquad
 m(x+y) \leq C_0\, m(x) m(y).
\end{equation}
Common examples are (for $\alpha_i\geq0$)~:
$$m_1(x)=1+|x_1|^{\alpha_1}+\ldots+|x_d|^{\alpha_d}
\qquad\text{and}\qquad m_2(x)=e^{\alpha|x|}.$$
The class of Peetre-type weights is stable by finite sums and products,
translations and orthogonal transforms.

\begin{lemma}\label{peetre-prop}
Let $m$ be a Peetre-type weight such that $m(x)\leq C \exp(c|x|)$
and~$T>0$. Then, there is a constant $C_T>0$ such that
\begin{equation}
\norme[L^\infty]{m (e^{t\Delta}a)} \leq C_T \norme[L^\infty]{m a}.
\end{equation}
\end{lemma}
\Proof
It is an elementary computation~:
\begin{align*}
m(x) \left|e^{t\Delta}a(x)\right| &\leq C_0 \left[(m g_t) \ast (m |a|)\right](x)\\[3pt]
&\leq C_0 \, (4\pi)^{-d/2} \norme[L^\infty]{ma} \int_{\R^d} m\big(\sqrt{t}\,y\big)
\: e^{-y^2/4}\,dy.
\end{align*}
The conclusion follows from the bound $m\big(\sqrt{t}\,y\big) \leq C \exp(c\, T\, |y|)$.

\endProof

\bigskip
As a converse to the previous result, the following property implies that
highly localized flows cannot decay at infinity in a really anisotropic way.

\begin{proposition}\label{corollary3}
Let $a\in L^\infty_{d+1+\varepsilon}$ be a divergence-free vector field
with $0<\varepsilon<1$, and  $u$ the corresponding solution of \NS in
${\cal C}_w\left([0,T);L^\infty_{d+1}\right)$.
For some $t>0$, let us assume that there exist an index $j\in\{1,\ldots,d\}$
and a subset $\Sigma\subset\mathbb{S}^{d-1}$ of positive measure such that
\begin{equation}\label{high-anisotropic}
\forall \sigma\in\Sigma, \qquad
\lim_{\substack{|x|\to+\infty\\ x\in\R\sigma}} |x|^{d+1} u_j(x,t) = 0.
\end{equation}
Then, there exists a constant $C>0$ such that
\begin{equation}
|u_k(x,t)|\le C(1+|x|)^{-(d+1+\varepsilon)}
\end{equation}
for all $k=1,\ldots,d$. Moreover, if \eqref{high-anisotropic} holds for
a finite time interval $t\in[T_0,T_1]$, then $C$ may be chosen
uniformly with respect to $t$.
\end{proposition}

\Proof
Our assumptions imply that the polynomial $P_j(x,t)$ identically vanishes.
Proposition~\ref{theorem2} then implies that 
all the other components of $P(x,t)$ also vanish.
Our statement is once again a consequence of~\eqref{profile}.

\endProof

\subsection{Application to the decay in a half-space domain}
\label{sub-half}

Our last application of Theorem~\ref{theorem1} is the study of the decay of solutions
of the Navier--Stokes equations in the half space
$$\R^d_+=\{(x',x_d)\;:\;x'\in\R^{d-1},\,x_d>0\}.$$
We set $u'=(u_1,\ldots,u_{d-1})$ and $x'=(x_1,\ldots,x_{d-1})$.
Let $\{e^{-tA'}\}_{t\ge0}$ be the semigroup generated by $-A'=\Delta$, 
in the case of the Neumann boundary conditions~:
\begin{equation}
\label{neumann}
\partial_d u'|_{\partial \R^d_+}=0, \qquad u_d|_{\partial \R^d_+}=0.
\end{equation}
where $\partial_d=\frac{\partial}{\partial x_d}$.
The integral formulation of the Navier--Stokes system in $\R^d_+$ is
\begin{equation}
\label{IEhalf}
u(t)=e^{-tA'}a-\int_0^te^{-(t-s)A'}\P\div(u\otimes u)(s)\,ds,
\end{equation}
with $\hbox{div\,}a=0$.
We refer to \cite{FujM02} for the construction of weak and strong
solutions to~\eqref{IEhalf}.
\bigskip

We have the following result:

\begin{proposition} \label{corollary2}
Assume that $a\in L^\infty_\vartheta(\R^d_+)$, with $\vartheta>\frac{(d+1)}{2}$.
Then there exist $T>0$ and a unique
strong solution $u\in {\cal C}_w\left([0,T);L^\infty_{\ubar\vartheta} (\R^d_+)\right)$
of  \eqref{IEhalf}.
Such a solution satisfies
\begin{equation}
\label{profile-half}
u(x,t)=e^{t\Delta}a(x)+H(x,t)
+{\cal O}_t\left(|x|^{-\min\{2\vartheta \,;\, d+2\}}\right),
\end{equation}
where $H=(H_1,\ldots,H_d)$ is homogeneous of degree $-(d+1)$  for all $t\in[0,T)$,
and such that:
\begin{eqnarray}
\label{boundQ-j}
|H_j(x,t)|&\le& C |x'|\cdot |x|^{-(d+2)}, \qquad (1\le j\le d-1)\\
\label{boundQ-d}
|H_d(x,t)|&\le& C |x_d|\cdot |x|^{-(d+2)}.
\end{eqnarray}
\end{proposition}

As an immediate consequence of Proposition~\ref{corollary2},
 we obtain the following anisotropic decay estimates
(assuming that $a$ is well localized).
\begin{equation*}
\begin{array}{llll}
u_d(x,t)=\phantom{\Bigl|\!} {\cal O}\bigl(|x'|^{-(d+2)}\bigr),& u'(x,t)={\cal O}\bigl(|x'|^{-(d+1)}\bigr) , \
	&& \mbox{when $|x'|\to+\infty$, $x_d$ fixed}\\
u_d(x,t)={\cal O}\bigl(|x_d|^{-(d+1)}\bigr), & u'(x,t)={\cal O}\bigl(|x_d|^{-(d+2)}\bigr) ,
&& \mbox{when $x_d\to+\infty$, $x'$ fixed}. \\
\end{array}
\end{equation*}
It is also worth noticing that Proposition~\ref{corollary3} is not violated as the
above decay holds only in a cylindrical region, and not in a conical one.

\Proof
This is immediate. Indeed, the study of \eqref{IEhalf} is reduced to that of \NS
in the following way. If $u$ solves \eqref{IEhalf}, then one can construct
a solution of \NS in the whole $\R^d$, setting $$\tilde
u_j(x_1,\ldots,x_{d-1},-x_d,t)=u_j(x_1,\ldots,x_{d-1},x_d,t)$$
for $j=1,\ldots,d-1$ and $\tilde u_d(x_1,\ldots,x_{d-1},-x_d,t)=-u_d(x_1,\ldots,x_{d-1},x_d,t)$
(see \cite{FujM02}).
Then under the assumptions of Proposition~\ref{corollary2} we can apply  \eqref{profile} to $\tilde u$.
But the integrals $\tilde \cc_{j,d}(t)\equiv \int_0^t\!\!\int_{\R^d}(\tilde u_j\tilde u_d)(x,s)\,dx\,ds$ vanish,
for $j\not=d$.
Hence, from \eqref{polynomial}
we see that $|P(x,t)|$ is bounded by a function $H(x,t)$ satisfying 
\eqref{boundQ-j}. 

\endProof

\subsection*{Conclusions}

Theorem~\ref{theorem1}
provides a quite complete answer to the spatial decay
problem of solutions of the free Navier--Stokes equations in the whole space,
at least for well localized data.
It would be interesting to know if some of the results of the present paper
can be adapted to flows in other domains.

\medskip
For example, in the half-space case, the Neumann boundary condition
considered in the previous section is not the most interesting one,
since it destroys the boundary layer effects.
The construction of the asymptotics as in Theorem~\ref{theorem1}, in the
case of Dirichlet boundary conditions, would require a careful analysis of 
Ukai's formula (or its more recent reformulations)
for the Stokes semigroup.

In the case of stationary flows, asymptotic profiles have been given,
\eg by F.~Haldi and P.~Wittwer \cite{HalW05}, \cite{Wit02}.
Their results model the wake flow beyond an obstacle. However, they
do not deal with the obstacle itself, but with a half-plane domain
and a technical boundary condition dictated by experimental knowledge.

\medskip
For the non-stationary equation \eqref{NS} in $\R^3\backslash\Omega$
with Dirichlet boundary conditions on~$\partial\Omega$,
it seems reasonable to expect
that anisotropic lower bound estimates for the decay of $u$ should hold,
when the net forces exerted by the fluid on the boundary, \ie
$$ \int_{\partial_\Omega}\bigl(T[u,p]\cdot \nu\bigr)(y,t)\,dS_y$$
(where $T_{j,k}[u,p]=\partial_ju_k+\partial_ku_j-\delta_{j,k}p$
and $T[u,p]=(T_{j,k}[u,p])_{j,k}$ is the stress tensor)
do not vanish. This last condition, which is motivated by
the results of Y.~Kozono  \cite{Koz98} and C.~He,  T.~Miyakawa \cite{HeM03}
on the $L^1$-summability of solutions,
would play, in the exterior domain case, a role equivalent to
the non-vanishing criterion given by Proposition~\ref{theorem2}.

%%%%%%%%%%%%%%%%%%%%%%%%%%%%%%%%%%%%%%%%%%%%%%%%%%%%%%%%%%%%%%%%%%
\bibliographystyle{amsplain}

\end{document}